\documentclass[12pt]{article}
\newtheorem{theorem}{Theorem}[section]
\newtheorem{lemma}{Lemma}[section]

\usepackage{epsfig}
\usepackage{color}
\usepackage{amsmath}
\usepackage{setspace}
\usepackage{array}
\usepackage{color}
\usepackage{amssymb}
\usepackage{amsbsy}
\usepackage{multirow}

\title{D-optimal Factorial Designs under Generalized Linear Models}

\author{Jie Yang$^{1}$ and Abhyuday Mandal$^{2}$\\
$^1$University of Illinois at Chicago and $^2$University of Georgia}

\begin{document}
\maketitle

\begin{quotation}
\begin{center} {\bf{ \it Abstract:}}
\end{center}
\bigskip\noindent
Generalized linear models (GLMs) have been used widely for modelling the mean response both for discrete and continuous random variables with an emphasis on categorical response. Recently Yang, Mandal and Majumdar (2013) considered full factorial and fractional factorial locally D-optimal designs for binary response and two-level experimental factors. In this paper, we extend their results to a general setup with response belonging to a single-parameter exponential family and for multi-level predictors.

\bigskip
\par

\par
\end{quotation}

\section{Introduction}

\bigskip\noindent Binary responses and count data are usually modelled using generalized linear models (GLMs). GLMs have been widely used for modelling the mean response both for discrete and continuous random variables with an emphasis on categorical response. Although the methods of analyzing data using these models have been discussed in depth in the literature (McCullagh and Nelder (1989), Dobson and Barnett (2008)), only a limited number of results are available for the choice of optimal design of experiments under GLMs (Khuri et al. (2006), Woods et al. (2006), Atkinson et al. (2007), Stufken and Yang (2012)). For optimal designs under GLMs, there are four different approaches proposed in the literature to handle the dependence of the design optimality criterion on the unknown parameters, (1) local optimality approach of Chernoff (1953) in which the parameters are replaced by assumed values; (2) Bayesian approach (Chaloner and Verdinelli (1995)) that considers a prior belief on unknown parameters; (3) maximin approach that maximizes the minimum efficiency over range of values of the unknown parameters (see Prozanto and Walter (1988) and Imhof (2001)); and (4) sequential approach where the estimates of the design parameters are updated in an iterative way (see Ford, Titterington and Kitsos (1989)). 

\bigskip\noindent In this paper, we will focus on local optimality and study D-optimal factorial designs under GLMs. Recently Yang, Mandal and Majumdar (YMM hereafter) (2012) obtained some theoretical results on locally D-optimal factorial designs for $2^2$ experiments with binary response. In a subsequent paper,  YMM (2013) considered D-optimal factorial designs with  $k$ two-level predictors and binary response. In this paper, we extend their results to a much more general setup with multi-level predictors and a response belonging to a single-parameter exponential family.

\bigskip\noindent Following YMM (2012, 2013), we consider locally D-optimal factorial designs. In this case, the Fisher's information matrix contains the unknown parameters, so does the D-optimality criterion. The optimality criterion can be written in terms of the variances, or measures of information, at each of the design points. Note that these variances depend on the parameters through the link function. The locally D-optimal designs can be quite different from the {\it uniform design} (that is, the experimental units are distributed evenly among the pre-specified design points) which is commonly used in practice, especially when some of the regression parameters stay significantly away from zero (see YMM (2013)). For illustration suppose we consider an experiment with two two-level factors, where the response is modeled by a Poisson regression, which leads to the linear predictor $\eta=\beta_0+\beta_1x_1+\beta_2x_2$. For example, an insurance company might be interested in conducting an experimental study to count the number of car break-ins in Detroit during some time period. It is a $2^2$ experiment with factors {\tt parking location} ($x_1 = +1$ for ``off-street" or $-1$ for ``inside garage") and whether or not the car is equipped with any kind of {\tt anti-theft device} ($x_2=+1$ for ``Yes" or $-1$ for ``No"). If we have solid reasons to believe that {\tt anti-theft device} is a more dominant factor than {\tt parking location}, and  based on a previous study, an initial choice of parameters $(\beta_0, \beta_1, \beta_2) = (1, 1, -2)$ is reasonable, then the uniform design  has only 78.7\% relative efficiency compared with the D-optimal design (see Example~5.3 for the definition of relative efficiency).

\bigskip\noindent This paper is organized as follows. In Section~2, we provide a preliminary setup for the generalized linear models under consideration. We discuss the characterization of locally D-optimal designs in Section~3. In Section~4 we discuss our search algorithms, both theoretically and numerically, for obtaining D-optimal approximate or exact designs. In Section~5, we illustrate our results with some real examples. We conclude with some remarks in Section~6.

\section{Preliminary Setup: Generalized Linear Models}

\bigskip\noindent In this paper, we consider experiments with a univariate response $Y$ which follows an exponential family with a single parameter $\theta$ in canonical form
$$\xi(y;\theta) = \exp\{y b(\theta) + c(\theta) + d(y)\}~,$$
where $b(\theta), c(\theta)$ and $d(y)$ are known functions. For examples,
\[
\xi(y; \theta) = \left\{\begin{array}{ll}
\exp\left\{y\log\frac{\theta}{1-\theta} + \log(1-\theta)\right\}, & \mbox{binary }Y\sim {\rm Bernoulli}(\theta)\\
\exp\left\{y\log\theta - \theta - \log y!\right\}, & \mbox{count response }Y\sim {\rm Poisson}(\theta)\\
\exp\left\{y\frac{-1}{\theta} - k\log\theta + \log\frac{y^{k-1}}{\Gamma(k)}\right\}, & Y\sim {\rm Gamma}(k, \theta), k>0\mbox{ fixed}\\
\exp\left\{y\frac{\theta}{\sigma^2} - \frac{\theta^2}{2\sigma^2} - \frac{y^2}{2\sigma^2} - \frac{1}{2}\log(2\pi\sigma^2)\right\}, & y\sim N(\theta, \sigma^2), \sigma^2 > 0\mbox{ fixed}.
\end{array}\right.
\]

\bigskip\noindent Consider independent random variables $Y_1, \ldots, Y_n$ with corresponding covariates or predictors ${\mathbf x}_1, \ldots, {\mathbf x}_n$, where ${\mathbf x}_i = (x_{i1}, \ldots, x_{id})'\in \mathbb{R}^d$. Under a generalized linear model, there exists a link function $g$ and parameters of interest $\boldsymbol{\beta} = (\beta_1, \ldots, \beta_d)'$, such that
$$E(Y_i) = \mu_i\mbox{ and } \eta_i = g(\mu_i)= {\mathbf x}_i^T\boldsymbol{\beta}.$$
It is known that, for example, see McCullagh and Nelder (1989) or Dobson and Barnett (2008),
\begin{eqnarray*}
E(Y_i) &=&  -\frac{c'(\theta_i)}{b'(\theta_i)}\ ,\\
{\rm var} (Y_i) &=& \frac{b''(\theta_i)c'(\theta_i) - c''(\theta_i)b'(\theta_i)}{b'(\theta_i)^3}\ ,
\end{eqnarray*}
provided $b(\theta)$ and $c(\theta)$ are twice differentiable. Since $Y_i$'s are independent, the $d\times d$ Fisher information matrix ${\mathbf F} = (F_{jk})$ can be obtained by
$$F_{jk} = E\left(\frac{\partial l}{\partial \beta_j}\frac{\partial l}{\partial \beta_k}\right) = \sum_{i=1}^n \frac{x_{ij}x_{ik}}{{\rm var}(Y_i)}\left(\frac{\partial \mu_i}{\partial \eta_i}\right)^2,$$
where $l=\log f(Y;\theta)$ could be written as a function of  $\boldsymbol{\beta}$ too.

\bigskip
Suppose there are only $m$ distinct predictor combinations ${\bf x}_1, \ldots, {\bf x}_m$
with numbers of replicates $n_1, \ldots, n_m$, respectively. Then
$$F_{jk} = \sum_{i=1}^m \frac{n_ix_{ij}x_{ik}}{{\rm var}(Y_i)}\left(\frac{\partial \mu_i}{\partial \eta_i}\right)^2= n\sum_{i=1}^m x_{ij}\frac{p_i}{{\rm var}(Y_i)}\left(\frac{\partial \mu_i}{\partial \eta_i}\right)^2x_{ik}
$$
where $p_i = n_i/n$, $i=1,\ldots, m$.
That is, the information matrix ${\mathbf F}$ can written as
$${\mathbf F}=nX^TWX$$
where $X=({\mathbf x}_1, \ldots, {\mathbf x}_m)^T$ is an $m\times d$ matrix, and $W={\rm diag}(p_1w_1,\ldots, p_mw_m)$ with
$$w_i = \frac{1}{{\rm var}(Y_i)}\left(\frac{\partial \mu_i}{\partial \eta_i}\right)^2.$$
For typical applications, ${\mathbf F}$ or $X^TWX$ is nonsingular.

\bigskip\noindent Suppose the link function $g$ is one-to-one and differentiable. Further assume that $\mu_i$ itself determines ${\rm var}(Y_i)$, that is, there exists a function $h$ such that ${\rm var}(Y_i) = h(\eta_i)$. Let $\nu =  \left(\left(g^{-1}\right)'\right)^2/h$. Then $w_i = \nu(\eta_i) = \nu\left({\mathbf x_i}'\boldsymbol\beta\right)$ for each $i$. We illustrate it using the four cases below.

\bigskip\noindent
{\bf Example 2.1: Binary response} Suppose $Y_i \sim $Bernoulli($\mu_i$). Then $h = g^{-1}(1-g^{-1})$.
For commonly used link function $g$,
\begin{equation*}\label{hfunctionfig}
\nu(\eta)=\left\{
\begin{array}{ll}
\frac{1}{2 + e^\eta + e^{-\eta}}\ , & \mbox{ for logit link; }\\
& \\
\frac{\phi(\eta)^2}{\Phi(\eta)[1-\Phi(\eta)]}\ , & \mbox{ for probit link; }\\
& \\
\left(\exp\{e^{\eta}\}-1\right)\left[\log\left(1-\exp\{-e^\eta\}\right)\right]^2\ , & \mbox{ for complementary log-log link; }\\
& \\
\frac{\exp\{2\eta-e^\eta\}}{1-\exp\{-e^\eta\}}\ , & \mbox{ for log-log link. }
\end{array}
\right.
\end{equation*}

\bigskip\noindent
{\bf Example 2.2: Poisson count} Suppose $Y_i \sim $Poisson($\mu_i$). For the canonical link function $g=\log$, $h=\nu=\exp$ and $w_i=\nu(\eta_i)=\exp\{\eta_i\}$ for each $i$.

\bigskip\noindent {\bf Example 2.3: Gamma response} Suppose $Y_i \sim {\rm Gamma}(k, \mu_i/k)$ with known $k>0$. For the canonical link function $g(\mu)=1/\mu$, $h(\eta) = 1/(k\eta^2)$, $\nu(\eta) = k/\eta^2$ and $w_i = \nu(\eta_i) = k/\eta_i^2$. A special case is $k=1$ which corresponds to the exponential distribution.

\bigskip\noindent {\bf Example 2.4: Normal response} Suppose $Y_i \sim N(\mu_i, \sigma^2)$ with known $\sigma^2 > 0$. For the canonical link function $g(\mu)\equiv\mu$, $h(\eta)\equiv \sigma^2$, $\nu(\eta) \equiv 1/\sigma^2$ and $w_i \equiv 1/\sigma^2$.

\section{Locally D-optimal Designs}

\bigskip\noindent In this paper, we consider experiments with an $m\times d$ design matrix $X = ({\bf x}_1, \ldots, {\bf x}_m)^T$ which consists of $m$ distinct experimental settings, that is, $m$ distinct combinations of values of $d$ covariates (for many applications, the first covariate is just the constant $1$). The response $Y_i$ follows a single-parameter exponential family under the generalized linear model $\eta_i = g(E(Y_i)) = {\bf x}_i^T\boldsymbol{\beta}$. For instance, in a $ 2^{3} $ experiment with binary response, $\eta =\beta _{0}+\beta _{1}x_{1}+\beta _{2}x_{2}+\beta _{3}x_{3}+\beta _{23}x_{2}x_{3}$ represents a model that includes all the main effects and the two-factor interaction of factors $2$ and $3$. The aim of the experiment is to obtain inferences about the parameter vector of factor effects $\boldsymbol\beta ;$ in the preceding example, $\boldsymbol\beta =\left( \beta _{0},\beta _{1},\beta _{2},\beta _{3},\beta _{23}\right) ^{\prime }$.

\bigskip\noindent In the framework of locally optimal designs, we assume that $w_i = \nu({\bf x}_i'\boldsymbol\beta)$, $i$ $=$ $1, \ldots, m$, determined by regression coefficients $\boldsymbol\beta$ and link function $g$, are known. For typical applications, $w_i>0$, $i=1,\ldots,m$. Then the design problem we consider here is to find the ``optimal'' allocation of $n$ experiment units into $(n_1, \ldots, n_m)'$ such that $n_1+\cdots + n_m=n$, known as {\it exact design}, or alternatively, find the ``optimal" proportion $p_i=n_i/n$, $i=1,\ldots, m$, known as {\it approximate design}. In this paper, we mainly focus on the approximate designs, that is, ${\mathbf p} = (p_1, \ldots, p_m)'$ such that $p_1 + \cdots + p_m = 1$. Since the maximum likelihood estimator of $\boldsymbol\beta $ has an asymptotic covariance matrix that is the inverse of $nX^{\prime }WX$, a {\it (locally) D-optimal design} is a ${\mathbf p}$ which maximizes $$f({\mathbf p}) = |X'WX|.$$ For example, for the $2^2$ experiment with main-effects model and binary response, $f({\mathbf p})=16(p_1 p_2 p_3 w_1 w_2 w_3 + p_1 p_2 p_4 w_1 w_2 w_4 + p_1 p_3 p_4 w_1 w_3 w_4 + p_2 p_3 p_4 w_2 w_3 w_4)$. See YMM (2012) for a detailed analysis of this problem. For general cases, the following lemma expresses that the objective function $f({\mathbf p}) = |X'WX|$ is an order-$d$ homogeneous polynomial of $p_1, \ldots, p_m$. This is useful in determining the optimal $p_i$'s. See Gonz\'{a}lez-D\'{a}vila, Dorta-Guerra and Ginebra (2007) and YMM (2013).

\begin{lemma}\label{XWX:lemma}
Let  $X[i_1,i_2,\ldots,i_d]$ be the $d\times d$ sub-matrix consisting of the $i_1\mbox{th}$, $\ldots$, $i_d\mbox{th}$ rows of the design matrix $X$. Then
$$f({\mathbf p})=|X'WX|=\sum_{1\leq i_1<\cdots<i_d\leq m}|X[i_1,\ldots,i_d]|^2 \cdot p_{i_1}w_{i_1}\cdots p_{i_d}w_{i_d}.$$
\end{lemma}
Based on Lemma~\ref{XWX:lemma}, YMM (2013) developed characterization theorems for D-optimal designs and
minimally supported (that is, the number of distinct supporting points is equal to the number of parameters)
D-optimal designs for two-level factors with $x_{ij} \in \{-1, 1\}$ (YMM (2013), Section~3.1). Note that the entries of $X$ in Lemma~\ref{XWX:lemma} can be any real number, which allows us to extend their results to multiple-level factors. Following YMM (2013) we define for each $i=1, \ldots, m$,
\begin{equation}\label{f_i(x)}
f_i(z) =f\left(\frac{1-z}{1-p_i}p_1,\ldots,\frac{1-z}{1-p_i}p_{i-1},z, \frac{1-z}{1-p_i}p_{i+1},\ldots, \frac{1-z}{1-p_i}p_{m}\right), \>\>\> 0\leq z\leq 1.
\end{equation}
Note that $f_i(z)$ is well defined for all ${\mathbf p}$ satisfying $f({\mathbf p}) > 0$. Lemma~7.1 and Lemma~7.2 of YMM (2013) could be applied for our case too. Thus we obtain theorems below to characterize locally D-optimal allocations and minimally supported D-optimal allocations for experiments with multiple-level factors and fairly general responses under GLMs. The proofs are presented in the Appendix.  Note that the design matrix $X$ here consists of $m\times d$ real-number entries, which generalizes the results in YMM (2013).

\begin{theorem}\label{theorem30}
Suppose $f\left({\mathbf p}\right)>0$. Then ${\mathbf p}$ is D-optimal if and only if for each $i=1, \ldots, m$, one of the two conditions below is satisfied:
\begin{itemize}
\item[(i)] $p_i=0$ and $f_i\left(\frac{1}{2}\right) \leq\frac{d+1}{2^{d}}f({\mathbf p})$;
\item[(ii)] $0 < p_i \leq \frac{1}{d}$ and $f_i(0) =\frac{1-p_i d}{(1-p_i)^{d}}f({\mathbf p})$.
\end{itemize}
\end{theorem}

\begin{theorem}\label{saturation:theorem}
Let ${\mathbf I} = \{i_1,\ldots,i_{d}\}\subset\{1,\ldots,m\}$ be an index set satisfying $|X[i_1,\ldots,$ $i_{d}]|\neq 0$. Then the saturated design satisfying $p_{i_1}=p_{i_2}=\cdots=p_{i_{d}}=\frac{1}{d}$ is D-optimal if and only if for each $i\notin {\mathbf I}$,
\[
\sum_{j\in {\mathbf I}} \frac{|X[\{i\}\cup {\mathbf I}\setminus\{j\}]|^2}{w_j} \leq \frac{|X[i_1,i_2,\ldots,i_{d}]|^2}{w_i}.
\]
\end{theorem}

\bigskip
\noindent Theorem~\ref{theorem30} and Theorem~\ref{saturation:theorem} are for $X$ with real-number entries and thus can be applied to multiple-level factors. We illustrate the results using the two examples below.

\bigskip\noindent {\bf Example 3.1: One two-level covariate and one three-level covariate:} Suppose the design matrix $X$ is given by
\[
X = \left(\begin{array}{rrr}
1 & 1 & 1\\
1 & 1 & 0\\
1 & 1 &-1\\
1 &-1 & 1\\
1 &-1 & 0\\
1 &-1 &-1
\end{array}\right)
\]
which consists of six supporting points: $(1,1), (1,0),$ $(1,-1),$ $(-1,1),$ $(-1,0),$ $(-1,$ $-1)$. Based on Theorem~\ref{theorem30}, an allocation ${\mathbf p}=(p_1,0,p_3,p_4,0,p_6)'$ is D-optimal if and only if $(p_1,p_3,p_4,p_6)'$ is D-optimal among the designs restricted on the four boundary points $(1,1)$, $(1,-1)$, $(-1,1)$, $(-1,-1)$, and
\begin{eqnarray*}
v_2 \geq \frac{p_1v_3p_4v_6 + p_1v_3v_4p_6 + v_1p_3p_4v_6 + v_1p_3v_4p_6 + 4v_1 v_3 p_4 p_6}{12(v_1p_3p_4p_6 + p_1v_3p_4p_6 + p_1 p_3 v_4 p_6 + p_1 p_3 p_4 v_6)},\\
v_5 \geq \frac{p_1v_3p_4v_6 + p_1v_3v_4p_6 + v_1p_3p_4v_6 + v_1p_3v_4p_6 + 4p_1 p_3 v_4 v_6}{12(v_1p_3p_4p_6 + p_1v_3p_4p_6 + p_1 p_3 v_4 p_6 + p_1 p_3 p_4 v_6)},
\end{eqnarray*}
where $v_i = 1/w_i$, $i=1,\ldots, 6$. In this situation, the second covariate reduces to a two-level factor. As the model has three parameters, one can further reduce the number of experimental settings. Theorem~\ref{saturation:theorem} states the necessary and sufficient conditions for a D-optimal design to have only three supporting points. For example, $p_1=p_2=p_4=1/3$ is D-optimal if and only if  $v_3 \geq v_1  +  4v_2$, $v_5 \geq v_1  +  v_2  +  v_4$ and $v_6 \geq 4v_1  +  4v_2  +  v_4$; $p_3=p_5=p_6=1/3$ is D-optimal if and only if  $v_1 \geq v_3  +  4v_5  +  4v_6$,  $v_2 \geq v_3  +  v_5  +  v_6$ and $v_4 \geq 4v_5  +  v_6$.

\bigskip\noindent {\bf Example 3.2: $2\times 3$ factorial design:} Suppose the design matrix
\begin{equation}\label{matrix:2by3}
X=\left(\begin{array}{rrrr}
 1 &  1 &  1 &  1\\
 1 &  1 &  0 & -2\\
 1 &  1 & -1 &  1\\
 1 & -1 &  1 &  1\\
 1 & -1 &  0 & -2\\
 1 & -1 & -1 &  1
\end{array}\right),
\end{equation}
where the four columns correspond to effects $I$, $A$ (two-level), $B_l$ (three-level, linear component), and $B_q$ (three-level, quadratic component). A minimally supported D-optimal design consists of four design points. 
Based on Theorem~\ref{saturation:theorem}, for example, $p_1=p_2=p_3=p_4=1/4$ is D-optimal if and only if $v_5 \geq v_1 + v_2 +v_4$ and $v_6 \geq v_1 + v_3 + v_4$; $p_2=p_3=p_4=p_5=1/4$ is D-optimal if and only if $v_1 \geq v_2 + v_4 + v_5$ and $v_6 \geq v_2 + v_3 + v_5$; $p_3=p_4=p_5=p_6=1/4$ is D-optimal if and only if $v_1 \geq v_3 + v_4 + v_6$ and $v_2 \geq v_3 + v_5 + v_6$.
{Note that $p_1=p_3=p_4=p_6=1/4$ is not D-optimal since the corresponding $|X'WX|=0$.}

\section{Searching for D-optimal Designs}

\bigskip\noindent YMM (2013) developed very efficient algorithms for searching locally D-optimal approximate designs or exact designs with binary response and two-level factors. Essentially the same algorithms could be used for maximizing $f({\mathbf p}) = |X'WX|$ for more general setup, as long as $X$ consists of $m$ distinct $d$-dimensional real-number vectors and the response belongs to a single-parameter exponential family.

\subsection{Lift-one algorithm for D-optimal approximate design}

\bigskip\noindent For finding ${\mathbf p}$ maximizing $f({\mathbf p})=|X'WX|$, the lift-one algorithm, as well as Lemma~7.1 and Lemma~7.2 of YMM (2013), could be applied to our cases, after limited modifications. Recall that $f_i(z)$ is defined in equation~(\ref{f_i(x)}).

\begin{lemma}\label{algo1lemma30} (YMM, 2013, Lemma~7.1)
Suppose $f\left({\mathbf p}\right)>0$. Then for $i=1, \ldots, m$,
\begin{eqnarray}\label{lem:f_i(x)}
f_i(z) = az(1-z)^{d-1}+b(1-z)^{d},
\end{eqnarray}
for some constants $a$ and $b$.
If $p_i>0$, $b=f_i(0)$,
$a=\frac{f\left({\mathbf p}\right)-b\left(1-p_i\right)^{d}}{p_i\left(1-p_i\right)^{d-1}}$;
otherwise,
$b=f\left({\mathbf p}\right)$,
$a=f_i\left(\frac{1}{2}\right)\cdot 2^{d}-b$.
Note that $a\geq 0$, $b\geq 0$, and $a + b>0$.
\end{lemma}

\begin{lemma}\label{algo1lemma31} (YMM, 2013, Lemma~7.2)
Let $l(z)=az(1-z)^{d-1}+b(1-z)^{d}$ with $0\leq z \leq 1$ and $a\geq 0, b\geq 0, a+b >0$.
If $a>bd$, then
$\max_z l(z)=\left(\frac{d-1}{a-b}\right)^{d-1}\left(\frac{a}{d}\right)^{d}\mbox{ at }
z=\frac{a-bd}{(a-b)d}\>\><1.$
Otherwise, $\max_z l(z)=b$ at $z=0$.
\end{lemma}

\bigskip\noindent {\bf Lift-one algorithm} (for multiple-level factors and single-parameter exponential family response)
\begin{itemize}
 \item[$1^\circ$] Start with arbitrary ${\mathbf p}_0=(p_1,\ldots,p_m)'$ satisfying $0<p_i<1$, $i=1,\ldots,m$. Compute $f\left({\mathbf p}_0\right)$.
 \item[$2^\circ$] Set up a random order of $i$ going through $\{1,2,\ldots,m\}$.
 \item[$3^\circ$] For each $i$, determine $f_i(z)$ as in expression~(\ref{lem:f_i(x)}). In this step, either $f_i(0)$ or $f_i\left(\frac{1}{2}\right)$ needs to be calculated according to Lemma~\ref{algo1lemma30}.
 \item[$4^\circ$] Define ${\mathbf p}_*^{(i)} =\left(\frac{1-z_*}{1-p_i}p_1,\ldots,\frac{1-z_*}{1-p_i}p_{i-1},z_*, \frac{1-z_*}{1-p_i}p_{i+1},\ldots, \frac{1-z_*}{1-p_i}p_{m}\right)'$, where $z_*$ maximizes $f_i(z)$ with $0\leq z\leq 1$ (see Lemma~\ref{algo1lemma31}). Note that $f({\mathbf p}_*^{(i)})=f_i(z_*)$.
 \item[$5^\circ$] Replace ${\mathbf p}_0$ with ${\mathbf p}_*^{(i)}$, $f\left({\mathbf p}_0\right)$ with $f({\mathbf p}_*^{(i)})$.
 \item[$6^\circ$] Repeat $2^\circ\sim 5^\circ$ until convergence, that is, $f({\mathbf p}_0)=f({\mathbf p}_*^{(i)})$ for each $i=1,\ldots, m$.
\end{itemize}

\bigskip\noindent
{\bf Convergence and performance of lift-one algorithm:} To guarantee the convergence, we may modify the lift-one algorithm as in YMM (2013, Section~3.3.1). A similar proof could be applied here to show that (1) if the lift-one algorithm or its modified version converges at ${\mathbf p}^*$, then ${\mathbf p}^*$ is D-optimal; (2) the modified lift-one algorithm is guaranteed to converge (YMM, 2013, Theorem~3.3).

\bigskip\noindent YMM (2013, Section~3.3.1) also compared the time cost of lift-one algorithm with commonly used nonlinear optimization algorithms including Nelder-Mead, quasi-Newton, conjugate gradient, and simulated annealing. Overall, lift-one algorithm could be 100 times faster than those algorithms. In this paper, we run more simulations to check how the computational time and number of non-zero $p_i$'s vary across different numbers of factors and ranges of parameters.

\begin{figure}[ht]\label{figure:timenonzero}
\caption{Performance of lift-one algorithm for $m=2^k$ and $d=k+1$}
\begin{center}
\includegraphics[scale=0.5]{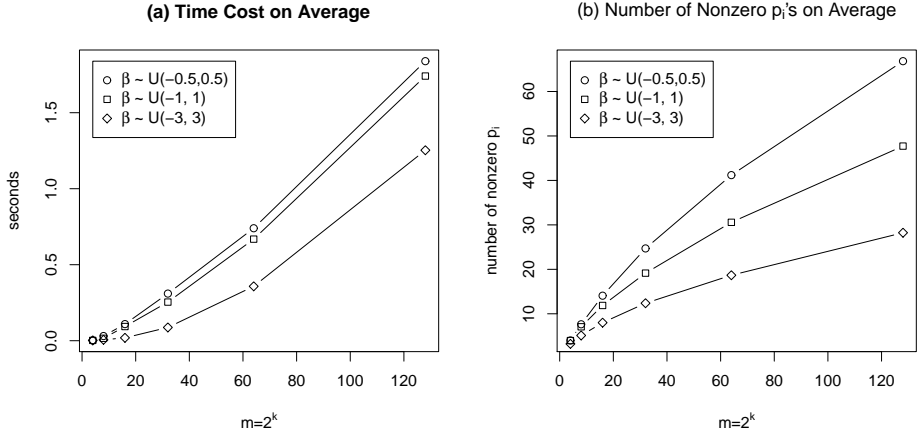}
\end{center}
\end{figure}

\bigskip\noindent Figure~1 shows the time cost and number of nonzero $p_i$'s on average based on 1000 simulated $\boldsymbol\beta$ for $2^k$ main-effects model with logit link. The time cost is based on a Windows Vista PC with Intel Core2 Duo CPU at 2.27GHz and 2GB memory. The relationship between the time cost and the number of supporting points $m=2^k$ is close to linear for moderate $k$. The time cost on average is only 1.25 secs, 1.74 secs, and 1.84 secs on average for $m=2^7=128$ with $\beta_i$'s follow iid  $U(-3,3)$, $U(-1,1)$, and $U(-0.5,0.5)$, respectively. In the meantime, the corresponding number of nonzero $p_i$'s on average are 28, 48, and 67 respectively.

\bigskip\noindent From a practitioner's point of view, it is often desirable to keep the number of supporting points of a design small. Due to Lemma~\ref{algo1lemma31}, the lift-one algorithm may force a $p_i$ to be exactly zero. It is an advantage of the lift-one algorithm over some commonly used nonlinear optimization algorithms which may keep an unimportant $p_i$ a tiny value but never let it to be zero. In many cases, it is hard to distinguish ``negligible'' from ``tiny''. As the range of $\beta_i$'s increases from $U(-0.5,0.5)$ to $U(-3, 3)$, the mean number of nonzero $p_i$'s reduces by one half due to more and more small $w_i$'s generated.

\subsection{Exchange algorithm for D-optimal exact designs}

\bigskip\noindent Given the total number of experimental units $n$, to find D-optimal integer allocation ${\mathbf n} = (n_1, \ldots, n_m)'$, YMM (2013) proposed another algorithm for two-level factors and binary response, called {\it exchange algorithm}, which adjusts $n_i$ and $n_j$ simultaneously for randomly chosen index pair $(i,j)$. The essentially same algorithm after some modifications could be used for a general design matrix $X$ consisting of real numbers.  The goal is to find the optimal ${\mathbf n}$ which maximizes $f({\mathbf n})=|X'W_nX|$, where $W_n = {\rm diag}\{n_1w_1,\ldots,n_mw_m\}$. To do this, we need a modified version of Lemma~7.5 and Lemma~7.4 in YMM (2013).

\begin{lemma}\label{algo3lemma2} (YMM, 2013, Lemma~7.5)
Fixing $1\leq i<j\leq m$, let
\begin{eqnarray}\label{lem:67fij}
f_{ij}(z)&=&\nonumber
f\left(n_1,\ldots,n_{i-1},z,n_{i+1},\ldots,n_{j-1},s-z,n_{j+1},\ldots,n_m\right)\\
&\stackrel{\triangle}{=}& Az(s-z)+Bz+C(s-z)+D,
\end{eqnarray}
where $s=n_i+n_j$.
Then
{\it (i)} $D>0 \Longrightarrow B>0\mbox{ and }C>0$;
{\it (ii)} $B>0\mbox{ or }C>0 \Longrightarrow A>0$;
{\it (iii)} $f({\mathbf n})>0\Longrightarrow A>0$;
{\it (iv)} $D = f(n_1,\ldots,$ $n_{i-1},0,n_{i+1},\ldots,$ $n_{j-1},0,n_{j+1},\ldots,$ $n_m)$.
{\it (v)} Suppose $s>0$, then
$A=\frac{2}{s^2}\left(2 f_{ij}\left(\frac{s}{2}\right)-f_{ij}(0)-f_{ij}(s)\right)$,
$B=\frac{1}{s}\left(f_{ij}(s)-D\right)$,
$C=\frac{1}{s}\left(f_{ij}(0)-D\right)$.
\end{lemma}

\begin{lemma}\label{algo3lemma1} (YMM, 2013, Lemma~7.4)
Let $q(z)=Az(s-z)+Bz+C(s-z)+D$ for real numbers $A > 0, B\geq 0, C\geq 0, D\geq 0$, and
integers $s>0, 0\leq z\leq s$.
Let $\Delta$ be the integer closest to $\frac{sA+B-C}{2A}$.
\begin{itemize}
\item[(i)] If $0\leq \Delta \leq s$, then
$\max_{0\leq z \leq s} q(z)=sC+D+(sA+B-C)\Delta-A\Delta^2
\mbox{ at }z=\Delta.$
\item[(ii)] If $\Delta < 0$, then $\max_{0\leq z\leq s} q(z)=sC+D$ at $z=0$.
\item[(iii)] If $\Delta > s$, then $\max_{0\leq z\leq s} q(z)=sB+D$ at $z=s$.
\end{itemize}
\end{lemma}

\bigskip
\noindent {\bf Exchange algorithm for D-optimal integer-valued allocations:}
\begin{itemize}
 \item[$1^\circ$] Start with initial design ${\mathbf n}=(n_1,\ldots,n_m)'$ such that $f({\mathbf n})=|X'W_nX|>0$.
 \item[$2^\circ$] Set up a random order of $(i,j)$ going through all pairs $$\{(1,2),(1,3),\ldots,(1,m), (2,3), \ldots,(m-1,m)\}.$$
 \item[$3^\circ$] For each $(i,j)$, let $s=n_i+n_j$. If $s=0$, let ${\mathbf n}^*_{ij}={\mathbf n}$. Otherwise, calculate $f_{ij}(z)$ as given in expression~(\ref{lem:67fij}). Then let
 $${\mathbf n}^*_{ij} =\left(n_1,\ldots,n_{i-1},z_*,n_{i+1},\ldots,n_{j-1},s-z_*,n_{j+1},\ldots,n_m\right)' $$
 where the integer $z_*$ maximizes $f_{ij}(z)$ with $0\leq z\leq s$ according to  Lemma~\ref{algo3lemma1}. Note that $f({\mathbf n}^*_{ij})=f_{ij}(z_*)\geq f({\mathbf n})>0$.
 \item[$4^\circ$] Replace ${\mathbf n}$ with ${\mathbf n}^*_{ij}$,  $f({\mathbf n})$ with $f({\mathbf n}^*_{ij})$.
 \item[$5^\circ$] Repeat $2^\circ\sim 4^\circ$ until convergence (no more increase in terms of $f({\mathbf n})$ by any pairwise adjustment).
\end{itemize}

The exchange algorithm usually converges in a few rounds. Nevertheless, it should be noted that the exchange algorithm for integer-valued allocations is not guaranteed to converge to the optimal ones. For further discussions on the convergence of the exchange algorithm for searching integer-valued solutions, see Section 3.3.2 of YMM (2013).

\section{Real Examples}

\bigskip \noindent In this section, we use some real examples to show how our results work.

\bigskip \noindent {\bf Example 5.1: Printed Circuit Board:}\hspace{0.2cm} Jeng, Joseph and Wu (2008) reported an experiment on printed circuit boards (PCBs), which has been modified to suit our purpose. This experiment was about inner layer manufacturing of PCBs. Several types of faults may occur during the manufacturing of PCBs, of which shorts and opens in the circuit are the major ones. In this case, we consider whether there is an open in the circuit as response, and two factors at two  and three levels respectively. The factors are whether preheating was done or not, and lamination temperature ($95^\circ$C, $105^\circ$C and $115^\circ$C). Table~5.2 is obtained from Table~2 of Jeng, Joseph and Wu (2008).

\bigskip
\begin{center}
Table 5.1: PCB Data (Jeng, Joseph and Wu, 2008)

\begin{tabular}{rrc}\hline
Preheat &   Temperature   &  \#success out of 480 replicates\\ \hline
1   &          1   &              120\\
1   &          2   &              16\\
1   &          3   &              25\\
2   &          1   &              50\\
2   &          2   &              51\\
2   &          3   &              22\\ \hline
\end{tabular}
\end{center}

\bigskip \noindent The design matrix is given by expression~(\ref{matrix:2by3}). Fitting a logistic regression model, we get
$\hat{\boldsymbol\beta} = ( -2.37,   0.154,   0.717,   0.113)^\prime.$ Assume that the true  ${\boldsymbol\beta}$ is about $(-2.5,0.15, 0.70,$ $0.10)^\prime$. Then the D-optimal design is ${\mathbf p}_o =(0.216, 0.186,$ $0.198,$ $0.206,$ $0.115,$ $0.080)'$. Applying the exchange algorithm for integer-valued allocations, we get the optimal allocation as $(621, 535, 569, 593, 331, 231)^\prime$ which is far away from being uniform, that is, each with 480 replicates.

\bigskip \noindent {\bf Example 5.2: Hard Disk Failures in a $2^2$ Experiment:}\hspace{0.2cm} Failure patterns of hard disk drives received considerable attention recently (Pinheiro, Weber and Barroso (2007)). Motivated by the data provided by Schroeder and Gibson (2007), we consider the following example, where the response can be modelled well by a Poisson regression. It is a $2^2$ experiment with factors as {\tt types of cluster} ($A$) and {\tt disk parameters} ($B$). The two levels of factor $A$ are HPC ($-1$) and IntServ ($+1$) where HPC represents clusters with organizations using supercomputers and IntServ represents large internet service providers. The two levels of factor $B$ are 10K rpm SCSI drives ($-1$) and 15K rpm SCSI drives ($+1$). After analyzing the data provided by Schroeder and Gibson (2007), one may consider the initial guess for ${\boldsymbol\beta}$ as $(5.5,-0.18, -0.22)^\prime$ with the design matrix
\[
X = \left(\begin{array}{rrr}
1 & 1 & 1\\
1 & 1 &-1\\
1 &-1 & 1\\
1 &-1 &-1
\end{array}\right).
\]
The D-optimal design for this choice of parameters is ${\mathbf p}_o =(0.18, 0.27, 0.26, 0.29)'$, which is close to the uniform design ${\mathbf p}_u=(0.25, 0.25, 0.25, 0.25)'$. Russell et al. (2009) characterized the locally D-optimal designs for Poisson regression models with the canonical log link. In their Remark 3, they mentioned a D-optimal design with the same supporting point as those of ours, for parameter ${\boldsymbol\beta}=(-0.91,0.04,-0.69)^\prime$ with optimal allocations ${\mathbf p}_o =(0.213,0.313,0.163,0.311)'$. Our algorithm yields exactly the same result, as expected.

\bigskip \noindent {\bf Example 5.3: Hard Disk Failures in a $2\times 3$ Experiment:}\hspace{0.2cm} In the context of Poisson regression, Anderson (2013) mentioned the number of hard disk failures at University of Illinois at Urbana Champaign  during a year. The administrators may be interested in finding the number of hard disk failures at a university during a year, in order to decide the amount of new purchases and adjust the budget for the next year. One may design a $2\times 3$ factorial experiment with the factors being {\tt types of computers} ($A$)  and {\tt operating system} ($B$). Both factors are categorical.  Desktop ($-1$) and Laptop ($+1$) are the two levels of factor $A$. Linux, Mac OS and Windows are the three levels of factor $B$. As suggested by Wu and Hamada (2009), one may choose the following two contrasts to represent the two degrees of freedom for the main effects of $B$:
\[
B_{01} = \left\{\begin{array}{rlll}
-1 & &0 &\\
1 &\mbox{for level} &1 &\mbox{of factor } B,\\
0 & &2 &
\end{array}\right. \nonumber\\
B_{02} = \left\{\begin{array}{rlll}
-1 & &0 &\\
0 &\mbox{for level} &1 &\mbox{of factor } B,\\
1 & &2 &
\end{array}\right. \nonumber
\]
where the levels $0$, $1$ and $2$ represent Linux, Mac OS and Windows respectively.
Then we have
\[
\mbox{Planning Matrix} =
\left(\begin{array}{cc}
\mbox{Desktop} & \mbox{Linux} \\
\mbox{Desktop} & \mbox{Mac OS} \\
\mbox{Desktop} & \mbox{Windows} \\
\mbox{Laptop} & \mbox{Linux} \\
\mbox{Laptop} & \mbox{Mac OS} \\
\mbox{Laptop} & \mbox{Windows}
\end{array}\right),
\]
and the design matrix as
\[
 X=\left(\begin{array}{rrrr}
 1 & -1 & -1 & -1\\
 1 & -1 &  1 &  0\\
 1 & -1 &  0 &  1\\
 1 &  1 & -1 & -1\\
 1 &  1 &  1 &  0\\
 1 &  1 &  0 &  1
\end{array}\right),
\]
where the four columns of $X$ correspond to effects $I$ (intercept), $A$ (two-level), $B_{01}$ (three-level), and $B_{02}$ (three-level) respectively. If one assumes that $\beta_0\sim U(-3,3)$,  $\beta_1\sim U(0,2)$, $\beta_2\sim U(0,1.5)$, and $\beta_3\sim U(0,3)$ independently, then the expected values of the $w_i$'s are $(0.24,  3.35,  9.18,  1.75, 24.76, 67.86)^\prime$. Following YMM (2013, Section~3.2), we calculate the EW D-optimal design ${\mathbf p}_o =(0,0,0.25,0.25,0.25,0.25)'$, which maximizes $|X'E(W)X|$. Conceptually, the EW D-optimality criterion is a surrogate to Bayes optimality (see YMM (2013) for more detailed discussion). In this case, it is interesting to note that the EW D-optimal design is a minimally supported one, which does not assign any observation to Desktops with Linux or Mac OS operating system. Compared with the EW D-optimal design ${\mathbf p}_o$, the relative efficiency of the uniform design ${\mathbf p}_u = (1/6, 1/6, \ldots, 1/6)'$ is 84\%. Note that the {\it relative efficiency} of ${\mathbf p}_u$ with respect to ${\mathbf p}_o$ is defined as $\left(f({\mathbf p}_u)/f({\mathbf p}_o)\right)^{1/d}$ in general.

\bigskip \noindent {\bf Example 5.4: Canadian Automobile Insurance Claims:}\hspace{0.2cm} Gamma regression can be used to model the total cost of car insurance claims given in Bailey and Simon (1960). They reported two studies for determining insurance rates using automobile industry data for policy years 1956$-$57. Motivated by their studies, a car insurance company may conduct an experiment to determine their own rates, with variables similar to those reported in Bailey and Simon (1960). Here one can consider a $2\times 4$ factorial example with factor levels given as follows:
\begin{center}
\begin{tabular}{lrl}
Factor & \multicolumn{2}{c}{Levels} \\
\hline
\multirow{2}{*}{Class} &		
  $+1$ &- pleasure\\
& $-1$ &- business use\\
\hline
\multirow{4}{*}{Merit} &
  3 &- licensed and accident free 3 or more years \\
& 2 &- licensed and accident free 2 years\\
& 1 &- licensed and accident free 1 year\\
& 0 &- all others \\
\hline
\end{tabular}
\end{center}
Let the design matrix be
\[
X=\left(\begin{array}{rrrrr}
 1 &  1 & 0 & 0 & 0   \\
 1 &  1 & 1 & 0 & 0   \\
 1 &  1 & 0 & 1 & 0   \\
 1 &  1 & 0 & 0 & 1   \\
 1 & -1 & 0 & 0 & 0   \\
 1 & -1 & 1 & 0 & 0   \\
 1 & -1 & 0 & 1 & 0   \\
 1 & -1 & 0 & 0 & 1
\end{array}\right),
\]
where the four columns correspond to effects $I$ (intercept), $A$ (two-level), $B_1$ (Merit1),  $B_2$ (Merit2), and $B_3$ (Merit3). After analyzing the data in the previous studies, one might consider $k = 1/55$ as described in Example~2.3. Motivated by their data, we set $\boldsymbol\beta = (-1, -0.75, -0.05, -0.25, -0.05)^\prime$ as the initial values of parameters. In this case, the optimal design is supported uniformly on the first and the last four rows of the design matrix mentioned above, that is, ${\mathbf p}_o =(.2,0,0,0,.2,.2,.2,.2)'$. The uniform design ${\mathbf p}_u = (1/8, 1/8, \ldots, 1/8)'$ has relative efficiency 82.7\% compared to the D-optimal ${\mathbf p}_o$.

\section{Discussion}

\bigskip \noindent In this paper, we extended the results of YMM (2013) to more general cases under the generalized linear model setup. Our framework allows the responses that belong to a single-parameter exponential family, which includes Binomial, Poisson, Gamma, exponential distributions as special cases. Our results also allow a fairly arbitrary set of design points, which could come from combinations of multiple-level factors, or grid points of continuous covariates.

\bigskip\noindent YMM (2013) also proposed EW D-optimal designs which is much easier to compute than the usual Bayesian D-optimal designs and more robust than uniform designs. The same concept and techniques can also be extended to a single-parameter exponential family response with arbitrary pre-specified design points. Their discussion on fractional factorial designs can also be extended here. This could be a topic of future research.

\section*{Acknowledgments}
The research of Abhyuday Mandal was in part supported by NSF Grant DMS-09-0573 and NSA Grant. The authors are thankful to Professor Dibyen Majumdar for many helpful discussions. They would also like to thank the reviewers for comments and suggestions that substantially improved the manuscript.

\section*{Appendix}

\noindent Note that we need Lemma~\ref{algo1lemma30} and Lemma~\ref{algo1lemma31} for the proof of Theorem~\ref{theorem30}. The proofs of Theorem~\ref{theorem30} and Theorem~\ref{saturation:theorem} here are similar to the ones in YMM (2013).

\noindent
{\bf Proof of Theorem~\ref{theorem30}:} Note that $f({\mathbf p}) > 0$ implies $0\leq p_i <1$ for each $i=1,\ldots, m$. Without any loss of generality, we assume $p_m >0$.

Define ${\mathbf p}_r=(p_1 ,\ldots,p_{m-1} )'$, and $f^{(r)}({\mathbf p_r})=f\left(p_1,\ldots, p_{m-1}, 1-\sum_{i=1}^{m-1}p_i\right)$. For $i=1, \ldots, m-1$, let $\boldsymbol\delta_i^{(r)}=(-p _1,\ldots,-p _{i-1},1-p _i,-p _{i+1},\ldots,-p _{m-1})'$. Then $f_i(z) = f^{(r)}({\mathbf p} _r + u\boldsymbol\delta_i^{(r)})$ with $u=\frac{z-p_i }{1-p_i }$. Since the determinant $|(\boldsymbol\delta_1^{(r)},\ldots,\boldsymbol\delta_{m-1}^{(r)})| = p_{m} \neq 0$,then $\boldsymbol\delta_1^{(r)}, \ldots, \boldsymbol\delta_{m-1}^{(r)}$ are linearly independent and thus form a basis to span the corresponding set of all feasible allocations $S_r = \{(x_1, \ldots, x_{m-1})'\ |\ \sum_{i=1}^{m-1} x_i \leq 1,\mbox{ and } x_i \geq 0, i=1, \ldots, m-1\}$ starting ${\mathbf p}_r$. Since $\log f^{(r)}$ is concave on the closed convex set $S_r$, then a solution maximizing $f^{(r)}$ exists and a local maximum of $f^{(r)}$ is also a global maximum on $S_r$. It can be verified that ${\mathbf p}_r$ maximizes $f^{(r)}$ on $S_r$ if and only if for each $i=1, \ldots, m-1$,
$\frac{\partial f^{(r)}({\mathbf p}_r +u\boldsymbol\delta_i^{(r)})}{\partial u}|_{u=0}
=0\mbox{ if }p_i >0\mbox{;}\>\>\>
\leq 0\mbox{ otherwise.}
$
That is, $f_i(z)$ attains its maximum at $z=p_i$, for each $i=1, \ldots, m-1$ (and thus for $i=m$).
Based on Lemma~\ref{algo1lemma30} and Lemma~\ref{algo1lemma31}, for $i=1, \ldots, m$, either
\begin{itemize}
\item[(i)] $p_i=0$ and $f_i\left(\frac{1}{2}\right)\cdot 2^d - f({\mathbf p}) \leq f({\mathbf p})\cdot  d$;\hspace{0.3cm} or
\item[(ii)] $p_i > 0$, $a > bd$, and $a- bd = p_i (a-b)d$, where $b = f_i(0)$,
and $a= \frac{f({\mathbf p}) - b(1-p_i)^d}{p_i(1-p_i)^{d-1}}$.
\end{itemize}
The conclusion can be obtained by simplifying those two cases.
\hfill{$\Box$}

\bigskip
{\bf Proof of Theorem~\ref{saturation:theorem}:}
Let ${\mathbf p}_I$ be the minimally supported design satisfying $p_{i_1}=p_{i_2}=\cdots=p_{i_d}=\frac{1}{d}$.
Note that if $|X[i_1,i_2,\ldots,i_d]| = 0$, ${\mathbf p}_I$ can not be D-optimal.
Suppose $|X[i_1,i_2,\ldots,i_d]|\neq 0$, ${\mathbf p}_I$ is D-optimal if and only if ${\mathbf p}_I$ satisfies the conditions of Theorem~\ref{theorem30}.
By Lemma~\ref{XWX:lemma},
$f({\mathbf p}_I)=d^{-d} |X[i_1,i_2,\ldots,i_d]|^2 w_{i_1}w_{i_2}$ $\cdots$ $w_{i_d}.$
\\
For $i\in {\mathbf I}$, $p_i=\frac{1}{d}$, $f_i(0)=0$. The condition~(ii) of Theorem~\ref{theorem30} is true. For $i\notin {\mathbf I}$, $p_i=0$, and
\begin{eqnarray*}
f_i\left(\frac{1}{2}\right) &=&(2d)^{-d} |X[i_1,\ldots,i_d]|^2 w_{i_1}\cdots w_{i_d}\\
                            &+&2^{-d}d^{-(d-1)} w_i\cdot w_{i_1}\cdots w_{i_d} \sum_{j\in {\mathbf I}} \frac{|X[\{i\}\cup {\mathbf I}\setminus\{j\}]|^2}{w_j}.
\end{eqnarray*}
Then $f_i\left(\frac{1}{2}\right) \leq f({\mathbf p}) \frac{d+1}{2^{d}}$ is equivalent to
$
\sum_{j\in {\mathbf I}} \frac{|X[\{i\}\cup {\mathbf I}\setminus\{j\}]|^2}{w_j} \leq \frac{|X[i_1,i_2,\ldots,i_d]|^2}{w_i}.
$
\hfill{$\Box$}

\section*{References}
\begin{enumerate}
\item Anderson, C. J. (2013). Poisson Regression or Regression of Counts (\& Rates), available at\\
{\tt http://courses.education.illinois.edu/EdPsy589/lectures/4glm3-ha-online.pdf}.

\item Atkinson, A. C., Donev, A. N. and Tobias, R. D. (2007). {\em Optimum Experimental Designs, with SAS}, Oxford University Press.

\item Bailey, R. and Simon, L. (1960). Two Studies in Automobile Insurance Rate Making. {\it ASTIN Bulletin I}, {\bf 4}, 192$-$217.

\item Chaloner, K. and Verdinelli, I. (1995). Bayesian experimental design: A review, {\em Statistical Science}, {\bf 10}, 273$-$304.

\item Chernoff, H. (1953). Locally optimal designs for estimating parameters, {\em Annals of Mathematical Statistics}, {\bf 24}, 586$-$602.

\item Dobson, A.J. and Barnett, A. (2008). \emph{An Introduction to Generalized Linear Models}, 3rd edition, Chapman and Hall/CRC.

\item Ford, I., Titterington, D. M. and Kitsos, C. P. (1989). Recent advances in nonlinear experimental design, {\it Technometrics}, {\bf 31}, 49$-$60.

\item Gonz\'{a}lez-D\'{a}vila, E., Dorta-Guerra, R. and Ginebra, J. (2007). On the information in two-level experiments, {\em Model Assisted Statistics and Applications}, {\bf 2}, 173$-$187.

\item Imhof, L. A. (2001). Maximin designs for exponential growth models and heteroscedastic polynomial models, {\em Annals of Statistics}, {\bf 29}, 561$-$576.

\item Jeng, S.-L., Joseph, V. R. and Wu, C.F.J. (2008). Modeling and analysis for failure amplification method, \emph{Journal of Quality Technology}, {\bf 40}, 128$-$139.

\item Khuri, A. I., Mukherjee, B., Sinha, B. K. and Ghosh, M. (2006). Design issues for generalized linear models: A review, {\em Statistical Science}, {\bf 21}, 376$-$399.

\item McCullagh, P. and Nelder, J. (1989). \emph{Generalized Linear Models}, 2nd edition, Chapman and Hall/CRC.

\item Pinheiro, E., Weber, W. and Barroso, L. A. (2007). Failure trends in a large disk drive population, {\it In Proceedings of the 5th USENIX Conference on File and Storage Technologies}, 17$-$28.

\item Pronzato, L. and Walter, E. (1988). Robust experiment design via maximin optimization, {\it Math. Biosci.}, {\bf 89}, 161$-$176.

\item Russell, K., Woods, D., Lewis, S. and Eccleston, J. (2009). D-optimal designs for Poisson regression models, {\it Statistica Sinica}, {\bf 19}, 721$-$730.

\item Schroeder, B. and Gibson, G. A. (2007). Understanding Disk Failure Rates: What Does an MTTF
of 1,000,000 Hours Mean to You? {\it ACM Transactions on Storage}, {\bf 3}, 8, 1$-$31.

\item Stufken, J. and Yang, M. (2012). Optimal designs for generalized linear models, In: \emph{Design and Analysis of Experiments, Volume 3: Special Designs and Applications}, K. Hinkelmann (ed.), Wiley, New~York.

\item Woods, D. C., Lewis, S. M., Eccleston, J. A. and Russell, K. G. (2006). Designs for generalized linear models with several variables and model uncertainty, {\em Technometrics}, {\bf 48}, 284$-$292.

\item Wu, C.F.J. and Hamada, M.S. (2009). \emph{Experiments: Planning, Analysis, and Optimization}, 2nd edition, Wiley.

\item Yang, J., Mandal, A. and Majumdar, D. (2012). Optimal designs for two-level factorial experiments with binary response, \emph{Statistica Sinica}, {\bf 22}, 885$-$907.

\item Yang, J., Mandal, A. and Majumdar, D. (2013). Optimal designs for $2^k$ factorial experiments with binary response, Technical Report, available at \\ {\tt http://arxiv.org/pdf/1109.5320v4.pdf}.
\end{enumerate}

\end{document}